\documentclass{amsart}
\usepackage{epsfig}
\usepackage{amssymb,mathrsfs}

\usepackage[frame,ps,matrix,arrow,curve,rotate,all,2cell,tips]{xy}
\input xy
\xyoption{all}

\usepackage[bookmarksopen, colorlinks, citecolor=blue, linkcolor=blue]{hyperref}

\usepackage[usenames,dvipsnames]{pstricks}

\usepackage{pst-grad} % For gradients
\usepackage{pst-plot} % For axes
\newtheorem{theorem}{Theorem}[section]
\theoremstyle{definition}

\theoremstyle{remark}
\newtheorem{remark}[theorem]{Remark}

\numberwithin{equation}{subsection}

\newcommand{\M}{\mathbb{M}}
\newcommand{\bbR}{\mathbb{R}}
\newcommand{\bbL}{\mathbb{L}}

\newcommand{\CC}{\mbox{$\mathcal C$}}

\def\VV{\mathbb{V}}

\def\Alg{\mathrm{Alg}}

\def\CC{\mathbb{C}}
\def\DD{\mathbb{D}}

\def\AA{\mathbb{A}}
\def\BB{\mathbb{B}}

\newtheorem{main}{Theorem}

\title{Substitudes, Bousfield localization, higher braided operads, and Baez-Dolan stabilization}

\author{David White}
\address{Department of Mathematics and Computer Science \\ Denison University
\\ Granville, OH 43023}
\email{david.white@denison.edu}

\begin{document}
\maketitle 
\begin{abstract}
This short note reports on joint work with Michael Batanin towards a general machine for proving Baez-Dolan Stabilization Theorems for various models of higher categories, based on substitudes, Bousfield localization, and homotopical Beck-Chevalley squares. I provide a road map to our recent papers, and include new results proving Baez-Dolan Stabilization Theorems for Tamsamani weak $n$-categories, higher Segal categories, Ara's $n$-quasi-categories, and cartesian models of Segal and complete Segal objects due to Bergner and Rezk. I also attempt to clarify the connection to higher braided operads, and our more general stabilization machinery.
\end{abstract}

\section{Introduction}

In 1995, Baez and Dolan introduced the \textit{stabilization hypothesis}, which loosely states that $k$-tuply monoidal weak $n$-categories are the same as $(k+1)$-tuply monoidal weak $n$-categories as long as $k\geq n+2$ \cite{BD}. Here \textit{$k$-tuply monoidal} signifies the additional structure you get on a weak $n$-category from reindexing from an $(n+k)$-category with one cell in each dimension $< k$. For example, if $\mathcal{C}$ is a 2-category with one object and one morphism, and we reindex two levels, then we obtain a 0-category (i.e., a set) with two commuting operations, corresponding to horizontal and vertical composition in the 2-cells of $\mathcal{C}$. By the Eckmann-Hilton argument, this yields the structure of a commutative monoid. Reindexing three levels, from a 3-category with only one cell in dimensions 0, 1, and 2, does not yield any additional structure on the resulting 0-category.

In \cite{batanin-white-CRM}, we sketched a proof of the stabilization hypothesis depending on the homotopy theory of $k$-operads (which encode $k$-tuply monoidal structure). We made good on this promise by proving the \cite[Theorem 14.2.1]{reedy}:

\begin{main}[Baez-Dolan Stabilization] \label{main:A}
Let  $0\le n$ and $\mathcal{M}$ an $n$-truncated monoidal combinatorial model category with cofibrant unit. Then 
  $i_!: B_k(\mathcal{M})\to B_{k+1}(\mathcal{M})$
  and
  $(j_k)_!: B_k(\mathcal{M})\to E_{\infty}(\mathcal{M})$
  are left Quillen equivalences for $k\ge n+2.$
\end{main}

Here $\mathcal{M}$ should be thought of as a model category of weak $n$-categories (e.g., Rezk's model via $\Theta_n$-spaces), $B_k(\mathcal{M})$ is the category of algebras over a $k$-operad $G_k$ (the cofibrant replacement of the terminal $k$-operad) encoding $k$-tuply monoidal weak $n$-categories, and $i$ and $j$ are comparison functors (based on suspension and symmetrization) between $k$-operads, $(k+1)$-operads, and symmetric operads, previously constructed by Batanin \cite{LocBat}. To say $\mathcal{M}$ is \textit{$n$-truncated} means its simplicial mapping spaces are $W_n$-local (defined below).

Such a result (but requiring a standard system of simplices on $\mathcal{M}$) had previously been proven by Batanin \cite{batanin-baez-dolan-via-semi}, but we deduce Theorem \ref{main:A} from a much stronger result \cite[Theorem 14.1.2]{reedy}, where $SO$ is the category of symmetric operads:

\begin{main}\label{main:B} Let $\mathcal{M}$ a combinatorial monoidal model category with cofibrant unit.
For $k\geq 3$ and $2\leq n+1\leq k$, the symmetrization functor
$sym_{k}^{}: Op^{W_n}_{k}(\mathcal{M}) \rightarrow SO(\mathcal{M})$
and the suspension functor 
$\Sigma_!:Op^{W_n}_{k}(\mathcal{M})\rightarrow  Op^{W_n}_{m}(\mathcal{M})$ (for $k<m\leq \infty$)
are {left Quillen equivalences}. Moreover, for $1\leq n\leq \infty$, the braided symmetrization functor $bsym_{2}^{}: Op^{W_n}_{2}(\mathcal{M}) \rightarrow BO(\mathcal{M})$
is a left Quillen equivalence with the category of braided operads.
\end{main} 

Here $Op^{W_n}_{k}(\mathcal{M})$ denotes the category of \textit{locally constant $k$-operads}, relative to the localizer $W_n$ that encodes $n$-types. As developed by Cisinski, a \textit{fundamental localizer} \cite[Definition 9.1.1]{reedy} is a class of functors between small categories that contains all identity functors, satisfies the two out of three property, is closed under retracts, contains functors $A \to 1$ where $1$ is the terminal category and $A$ is a category with terminal objects, and such that, if $u/c: A/c \to B/c$ is in $W$ for each object $c\in C$ (where $u$ is a morphism in $Cat/C$) then $u$ is in $W.$ 

The localizer $W_n$ is the smallest localizer containing the unique functor from the $(n+1)$-sphere (viewed as a category) to the terminal category. That minimal fundamental localizers such as $W_n$ exist is a theorem of Cisinski.  We recall that a category $A$ is said to be \textit{$W$-aspherical} if the unique functor from $A$ to $1$ is in $W$.

\section{Substitudes and left Bousfield localization}

To study the homotopy theory of $Op^{W_n}_{k}(\mathcal{M})$, we encode categories of $k$-operads as algebras over substitudes. A substitude \cite[Definition 5.1.1]{reedy} is equivalent to the data of a colored operad $P$ with a category $A$ of unary operations. This means one can encode structures with substitudes that cannot be encoded with colored operads. We use techniques from \cite{batanin-berger} and \cite{white-yau1} to transfer model structures from presheaf categories $[A,\mathcal{M}]$ to categories algebras over what we call $\Sigma$-free tame unary substitudes with faithful unit, a class that includes categories of $k$-operads. Notably, we prove a transfer theorem more general in two ways than those that have appeared previously. First, it works for substitudes rather than only for colored operads \cite[Theorem 8.1.7]{reedy}. Secondly, it works when the base, $\VV$, is only a semi-model category \cite[Theorem 2.2.1]{reedy}.

We generalize work of Cisinski \cite{cis} to prove the existence of left Bousfield localizations $[A,\mathcal{M}]^W$ for any proper fundamental localizer $W$, with respect to the projective, injective, or Reedy model structure on presheaves \cite[Theorem 9.3.5]{reedy}. In these local model structures, local objects $F: A\to \mathcal{M}$ are \textit{$W$-locally constant presheaves}, i.e., for any $W$-aspherical category $A'$, and any functor $u:A'\to A$, the induced functor $u^*(F): A'\to \mathcal{M}$ is isomorphic to a constant presheaf in $Ho[A',\mathcal{M}]$.  The local equivalences are morphisms $u: A\to B$ inducing right Quillen equivalences on categories of locally constant presheaves. 

$W_\infty$ is the minimal fundamental localizer making categories with terminal objects $W$-aspherical. Equivalently, $W_\infty$ is the class of functors whose nerve is a weak equivalence. $W_\infty$-locally constant functors $F$ are those taking all morphisms $f$ in $A$ to weak equivalences. This is analogous to \cite{chorny-white} where the local objects are the homotopy functors (i.e., those preserving weak equivalences). If $\mathcal{M}$ is $n$-truncated then $[A,\mathcal{M}]^{W_r} \to [A,\mathcal{M}]^{W_\infty}$ is a Quillen equivalence for all $r \geq n+1$.

To get from $Op_{k}(\mathcal{M})$ to $Op^{W_n}_{k}(\mathcal{M})$, we must left Bousfield localize. Unfortunately, categories of algebras over substitudes are often not left proper. To remedy this, we develop a theory of left Bousfield localization that does not require left properness, and results in a semi-model structure. A semi-model category \cite[Definition 2.1.1]{reedy} has three classes of morphisms that satisfy all of the model category axioms except that we only know that trivial cofibrations \textit{with cofibrant domain} lift against fibrations, and that morphisms \textit{with cofibrant domain} admit factorizations into trivial cofibrations followed by fibrations. Because semi-model categories admit cofibrant replacement, and because the subcategory of cofibrant objects behaves exactly like a model category, every result about model categories has a semi-model categorical analogue, and semi-model categories are equally useful in practice. We state our localization theorem \cite[Theorem A]{bous-loc-semi}:

\begin{main}[Bousfield localization without left properness] \label{bous-loc-semi}
Suppose that $\mathcal{M}$ is a combinatorial semi-model category whose generating cofibrations have cofibrant domain, and $\mathcal{C}$ is a set of morphisms of $\mathcal{M}$. Then {there is a semi-model structure $L_{\mathcal{C}}(\mathcal{M})$ on $\mathcal{M}$}, whose weak equivalences are the $\mathcal{C}$-local equivalences, whose cofibrations are the same as $\mathcal{M}$, and whose fibrant objects are the $\mathcal{C}$-local objects. Furthermore, $L_{\mathcal{C}}(\mathcal{M})$ {satisfies the universal property} that, for any any left Quillen functor of semi-model categories $F:\mathcal{M}\to \mathcal{N}$ taking $\mathcal{C}$ into the weak equivalences of $\mathcal{N}$, then $F$ is a left Quillen functor when viewed as $F:L_{\mathcal{C}}(\mathcal{M})\to \mathcal{N}$.
\end{main}

This theorem is of independent interest for a host of applications, detailed in \cite{bous-loc-semi}, as lack of left properness has bedeviled researchers seeking to left Bousfield localize for years. Examples in \cite{bous-loc-semi}, show that sometimes the classes of morphisms above do not satisfy the model category axioms, so only a semi-model structure is possible. Examples of semi-model structures abound \cite{batanin-baez-dolan-via-semi, batanin-white-CRM, batanin-white-eilenberg-moore, gutierrez-white-equivariant, hovey-white, white-thesis, white-user, white-loc, white-yau4, white-yau6, white-yau1, white-yau7, white-yau3, white-yau-5, white-yau-coloc}.

There are two ways to get from $[A,\mathcal{M}]$ to $Op^{W_n}_{k}(\mathcal{M})$. One can either localize first, then lift the resulting model structure (as in \cite{white-comm-mon, white-loc}), or one can lift first (using the transfer theorem) and then attempt to localize. As proven in \cite[Theorem 5.6]{batanin-white-eilenberg-moore}, these two approaches are equivalent (when both work).

In addition to these localization results, to prove Theorem \ref{main:B} we develop a theory of homotopical Beck-Chevalley squares \cite[Theorem 4.2.2]{reedy} to lift Quillen equivalences of presheaf categories to Quillen equivalences of algebras over substitudes. This vastly generalizing previous work on such problems (e.g., \cite{white-yau3}).

A square of right adjoints:

 \begin{equation*}\label{BScondition}\xymatrix@R=0.7em@C=0.7em{
\AA %\ar@<2.5pt>[r]^{\psi_!} 
%\ar@2{->}[rd]^b
%{\ar@{}[dr]|(.7){\Searrow}}
\ar@<0pt>[ddd]_{\beta^*}
& & &
\BB \ar@<0pt>[lll]_{\psi^*}
 \ar@<0pt>[ddd]^{\alpha^*} \\
  &\ar@2{->}[rd]^b &  & \\ 
  & & & \\
\CC %\ar@<2.5pt>[r]^{\phi_!}
% \ar@<2.5pt>[u]^{\beta_!}  
& & &
\DD \ar@<0pt>[lll]_{\phi^*} 
%\ar@<2.5pt>[u]^{\alpha_!}
}
\end{equation*} 

is called {\em Beck-Chevalley} if the natural transformation
\begin{equation*}\label{BC transformation} 
\mathbf{bc}:  \phi_!  \beta^* \to  \alpha^*  \psi_! 
\end{equation*}
is an isomorphism. This implies that if $(\phi_!,\phi^*)$ is an adjoint equivalence and $\beta^*,\alpha^*$ reflect isomorphims then $(\phi_!,\phi^*)$ is adjoint equivalence.

The homotopical version of this machinery says that the square above is {\em homotopy Beck-Chevalley} if $   \bbL \phi_! \bbR \beta^* (-)\to \bbR \alpha^* \bbL \psi_! (-)$ is an isomorphism in Ho$(\DD)$. This occurs if $\alpha^*$ preserves weak equivalences and $\beta^*$ preserves cofibrant objects. This implies that if $(\phi_!,\phi^*)$ is a Quillen equivalence and $\beta^*,\alpha^*$ reflect weak equivalences between fibrant objects, then $(\phi_!,\phi^*)$ is a Quillen equivalence.

We use this result to lift Quillen equivalences from presheaf categories $[A,\VV]_{proj}^W \leftrightarrows [B,\VV]_{proj}^W$ to algebras over substitudes $\Alg_P^W(\VV)\leftrightarrows \Alg_Q^W(\VV)$. Specifically, if a given a morphism of substitudes $(f,g):(P,A)\to (Q,B)$ induces a homotopical Beck-Chevalley square, we see that if $(g_!,g^*)$ is a Quillen equivalence, then so is $(f_!,f^*)$. A crucial ingredient in the proof of Theorem \ref{main:B} is that the morphisms comparing categories of $k$-operads, $(k+1)$-operads, and symmetric operads, do induce homotopical Beck-Chevalley squares, both before and after localization \cite[Proposition 11.3.2, Proposition 14.1.1]{reedy}.

\section{Higher Braided Operads}

Locally constant $k$-operads are a model for higher braided operads \cite{SymBat, LocBat}, and the category of unary operations $Q_k^{op}$ has $Q_k \cong \coprod Q_k(m)$ such that the nerve of $Q_k(m)$ is homotopy equivalent to the unordered configuration space of points in $\mathbb{R}^k$. An analysis of this homotopy type \cite[Theorem 11.1.7]{reedy} is the last ingredient in the proof of Theorems \ref{main:A} and \ref{main:B}, and the reason for the inequalities featuring $k$ and $n$. We also lift various equivalences of homotopy categories in this setting (known since \cite{LocBat}) to Quillen equivalences \cite[Proposition 12.2.1]{reedy}. Another consequence of Theorem \ref{main:B} is:

\begin{main}[Stabilization for Higher Braided Operads] \label{main:D}
If $\mathcal{M}$ is a $n$-truncated, combinatorial, monoidal model category with cofibrant unit, and $n\geq 0$ and {$3\leq n+2\leq k\leq \infty$}, then the symmetrization functor
$sym_{k}^{}: Op^{W_\infty}_{k}(\mathcal{M}) \rightarrow SO(\mathcal{M})$
and the suspension functor 
$\Sigma_!:Op^{W_\infty}_{k}(\mathcal{M})\rightarrow  Op^{W_\infty}_{m}(\mathcal{M})$ (for $k<m\leq \infty$)
are {left Quillen equivalences}.  Moreover, for $1\leq n\leq \infty$, $bsym_{2}^{}: Op^{W_\infty}_{2}(\mathcal{M}) \rightarrow BO(\mathcal{M})$
is a left Quillen equivalence.
\end{main}

This is proven in \cite[Corollary 14.1.3]{reedy}.

\section{Baez-Dolan Stabilization Theorems}

Finally, we obtain a stabilization result for $(n+m,n)$-categories, rather than just weak $n$-categories, as a consequence of the stronger results listed above. We state this first for Rezk's model of $(n+m,n)$-categories (where $Sp_m$ models $m$-types).

\begin{theorem} 
The suspension functor  induces a {left Quillen equivalence} 
$$ i_!: B_k(\Theta_n Sp_m) \to  B_{k+1}(\Theta_n Sp_m)$$  for  {$k\ge m+n+2$} and, hence, an  equivalence between homotopy categories of Rezk's $k$-tuply monoidal $(n+m,n)$-categories and Rezk's $(k+1)$-tuply monoidal $(n+m,n)$-categories. 
\end{theorem}

This is proven as \cite[Corollary 14.2.3]{reedy}, using that Rezk's $\Theta_n Sp_m$ is a $(n+m)$-truncated, combinatorial, monoidal model category with cofibrant unit. The same hypotheses apply to other models of higher categories, including Tamsamani weak $n$-categories, higher Segal categories, $n$-quasi-categories, and models of Bergner and Rezk for Segal and complete Segal objects in $\Theta_{n-1}$ spaces. These results are new, though we plan to add them to \cite{reedy}. 

\begin{theorem} \label{thm:tamsamani}
Let $\M$ be a combinatorial, monoidal model category with cofibrant unit. Then Simpson's categories $PC^n(\M)$ \cite[Theorem 19.3.2]{Simpson} are too, and hence satisfy our Stabilisation Theorem \ref{main:B}. If $\M$ is furthermore left proper, then Simpson's localization $Seg^n(\M)$, whose fibrant objects satisfy a Segal condition, satisfies our Stabilisation Theorem \ref{main:B}, and its $m$-truncation $\tau_m Seg^n(\M)$ satisfies our Theorem \ref{main:A}. In particular, the Baez-Dolan stabilization hypothesis is true for Tamsamani weak $n$-categories (corresponding to $\M = Set$, the trivial model structure) and higher Segal categories (corresponding to the Kan-Quillen model $\M = sSet$).
\end{theorem}

\begin{remark}
Simpson also proved a Baez-Dolan stabilization result \cite[Theorem 23.0.3]{Simpson}, but did not model $k$-tuply monoidal weak $n$-categories as algebras over a $k$-operad. Instead, he modeled them as $(k-1)$-connected weak $(n+k)$-categories. We conjecture that these two approaches are Quillen equivalent. Furthermore, we conjecture that `left proper' could be dropped above, using Theorem \ref{bous-loc-semi} to produce semi-model structures for the localizations, and then proving versions of Theorems \ref{main:A} and \ref{main:B} that only require a semi-model structure to begin. Lastly, it should be mentioned that Simpson requires $\M$ to be tractable, left proper, and cartesian (hence, to have cofibrant unit \cite[Definition 7.7.1]{Simpson}). Implicit in Theorem \ref{thm:tamsamani} is an extension of Simpson's approach to the realm of combinatorial, monoidal model categories. When $\M$ is monoidal but not cartesian, Simpson's connection to enrichments is lost, but one still has model categories $PC^n(\M)$ and $Seg^n(\M)$, and our Baez-Dolan stabilization result.
\end{remark}

We turn now to Ara's model category of $n$-quasi-categories, as presheaves over $\Theta_n$ that are Quillen equivalent to Rezk's model \cite[Theorem 8.4]{Ara}.

\begin{theorem}
Ara's model category of $n$-quasi-categories \cite[Theorem 2.2, Corollary 8.5]{Ara} satisfies the conditions of Theorem \ref{main:B} and its truncation $\tau_m nQcat$ satisfies the conditions of Theorems \ref{main:A} and \ref{main:D}. Hence, $n$-quasi-categories satisfy Baez-Dolan stabilization.
\end{theorem}

Lastly, we turn to two models introduced by Bergner and Rezk:
\begin{enumerate}
\item $\Theta_n Sp$-Segal categories, a combinatorial, cartesian model structure on functors $\Delta^{op}\to \Theta_n Sp$ whose fibrant objects satisfy a Segal condition \cite[Theorem 5.2]{BergnerRezk}. 
\item A combinatorial, cartesian model structure with all objects cofibrant, whose fibrant objects satisfy a subset of the conditions required of complete Segal spaces \cite[Proposition 5.9]{BergnerRezk}.
\end{enumerate}
Furthermore, they prove these two are equivalent to each other and to Rezk's $\Theta_n Sp$ \cite[Theorem 6.14, Corollary 7.1, Theorem 9.6]{BergnerRezk}. Lastly, both can be left Bousfield localized to make them $m$-truncated.

\begin{theorem}
Both of the Bergner-Rezk model structures above satisfy the conditions of Theorem \ref{main:B}, and their truncations satisfy the conditions of Theorems \ref{main:A} and \ref{main:D}. Hence, the Baez-Dolan stabilization hypothesis is true for these models.
\end{theorem}

These applications demonstrate the power of Theorems \ref{main:A}, \ref{main:B}, and \ref{main:D}: the conditions are satisfied by all models of higher categories that we are aware of that possess a monoidal model structure. There are other model of higher categories, including $n$-relative categories, $n$-fold Segal spaces, and simplicial categories, which may or may not possess a good monoidal product. If any of them is endowed with a monoidal model structure in the future, we anticipate that our methods will prove the Baez-Dolan stabilisation hypothesis for that model. Furthermore, any model that is homotopically equivalent to a model where we have proven the Baez-Dolan stabilization hypothesis will automatically satisfy stabilization on the homotopy category level.

\noindent {\bf Acknowledgments.} 
I am grateful to my co-author, Michael Batanin, for introducing me to the world of higher structures. I also thank the Oberwolfach Research Institute for Mathematics for hosting us for a week in September of 2021, which led to these new results and to this document. Lastly, 
I thank the organizers of that conference--Andrey Lazarev, Muriel Livernet, Michael Batanin, and Martin Markl--and 
I thank Alexander Campbell and Rhiannon Griffiths for questions that led to the new applications of our methods to Tamsamani weak $n$-categories and Ara's $n$-quasi-categories.

\end{document}